\documentstyle[11pt,twoside]{article}
\include{graphicx}
\title{A summation formula for a ${}_3F_2(1)$ hypergeometric series}
\author{\sc R. B.\ Paris \\
{\em Division of Computing and Mathematics}, \\
{\em Abertay University, Dundee DD1 1HG, UK}
}
\begin{document}
\def\f#1#2{\mbox{${\textstyle \frac{#1}{#2}}$}}
\def\dfrac#1#2{\displaystyle{\frac{#1}{#2}}}
\def\boldal{\mbox{\boldmath $\alpha$}}
{\newcommand{\Sgoth}{S\;\!\!\!\!\!/}
\newcommand{\bee}{\begin{equation}}
\newcommand{\ee}{\end{equation}}
\newcommand{\lam}{\lambda}
\newcommand{\ka}{\kappa}
\newcommand{\al}{\alpha}
\newcommand{\fr}{\frac{1}{2}}
\newcommand{\fs}{\f{1}{2}}
\newcommand{\g}{\Gamma}
\newcommand{\br}{\biggr}
\newcommand{\bl}{\biggl}
\newcommand{\ra}{\rightarrow}
\newcommand{\mbint}{\frac{1}{2\pi i}\int_{c-\infty i}^{c+\infty i}}
\newcommand{\mbcint}{\frac{1}{2\pi i}\int_C}
\newcommand{\mboint}{\frac{1}{2\pi i}\int_{-\infty i}^{\infty i}}
\newcommand{\gtwid}{\raisebox{-.8ex}{\mbox{$\stackrel{\textstyle >}{\sim}$}}}
\newcommand{\ltwid}{\raisebox{-.8ex}{\mbox{$\stackrel{\textstyle <}{\sim}$}}}
\renewcommand{\topfraction}{0.9}
\renewcommand{\bottomfraction}{0.9}
\renewcommand{\textfraction}{0.05}
\newcommand{\mcol}{\multicolumn}
\date{}
\maketitle
\pagestyle{myheadings}
\markboth{\hfill \sc R. B.\ Paris  \hfill}
{\hfill \sc  A hypergeometric sum\hfill}
\begin{abstract}
A summation formula is derived for the hypergeometric series of unit argument
${}_3F_2(1,1,c;d,n+2;1)$, where $n=0, 1, 2, \ldots$ and $\Re (d-c+n)>0$.
\vspace{0.4cm}

\vspace{0.3cm}

\end{abstract}

\vspace{0.3cm}

\noindent $\,$\hrulefill $\,$

\vspace{0.2cm}

\begin{center}
{\bf 1. \  Introduction}
\end{center}
\setcounter{section}{1}
\setcounter{equation}{0}
\renewcommand{\theequation}{\arabic{section}.\arabic{equation}}
In this note we give two proofs of the following summation theorem for a particular type of ${}_3F_2(1)$ hypergeometric series (also known as Clausen's series):
\newtheorem{theorem}{Theorem}
\begin{theorem}$\!\!\!.$\ Let $n=0, 1, 2, \ldots$ and $c$, $d$ be complex constants satisfying $\Re (d-c+n)>0$. Then
\[{}_3F_2\bl(\begin{array}{c}1,1,c\\d,n+2\end{array}\!\!;1\br)=\frac{(n+1) \g(d)}{(1-c)_{n+1}}\bl\{\frac{(d-c)_n}{\g(d-1)}\,[\psi(d-c+n)-\psi(d-1)]\]
\bee\label{e11}
-\sum_{k=1}^n \bl(\!\!\begin{array}{c}n\\k\end{array}\!\!\br)\bl(\!\!\begin{array}{c}n\!+\!1\!-\!c\\k\end{array}\!\!\br)
\frac{D_k(n)}{\g(d\!-\!n\!-\!1\!+\!k)}\br\},
\ee
where
\bee\label{e12}
D_k(n):=k!\{\psi(n+1)-\psi(n+1-k)\},
\ee
$(a)_k=\g(a+k)/\g(a)$ is Pochhammer's symbol and $\psi(a)=\g'(a)/\g(a)$ is the psi function.
\end{theorem}
The case $c=d$ (with $n\geq 1$) can be excluded from consideration since the resulting series contracts to a Gauss hypergeometric series which can be summed by the well-known Gauss summation theorem.
\vspace{0.6cm}

\begin{center}
{\bf 2. \ Proof 1}
\end{center}
\setcounter{section}{2}
\setcounter{equation}{0}
\renewcommand{\theequation}{\arabic{section}.\arabic{equation}}
In Miller and Paris \cite[(1.7)]{MP} the following summation theorem was established:
\[{}_3F_2\bl(\begin{array}{c}a,c,m\\d,m+p\end{array}\!\!;1\br)=(p)_m\sum_{k=0}^{p-1} (-)^k (m)_k \bl(\!\!\begin{array}{c}p-1\\k\end{array}\!\!\br)\,\frac{(1-d)_{k+m}}{(1-a)_{k+m}(1-c)_{k+m}}\]
\[+(m)_p \frac{\g(d) \g(d-a-c)}{\g(d-a) \g(d-c)} \sum_{k=0}^{m-1}(-)^k (p)_k \bl(\!\!\begin{array}{c}m-1\\k\end{array}\!\!\br)\,\frac{(d-a-c)_{k+p}}{(1-a)_{k+p} (1-c)_{k+p}}\]
for positive integers $m$ and $p$. In this formula we set $m=1$ and apply a limiting process to deal with the case $a\to 1$.
We therefore obtain
\[{}_3F_2\bl(\begin{array}{c}a,c,1\\d,p+1\end{array}\!\!;1\br)=\frac{p}{1-a}\sum_{k=0}^{p-1}(-)^k k!\bl(\!\!\begin{array}{c}p-1\\k\end{array}\!\!\br)\,\frac{(1-d)_{k+1}}{(1-c)_{k+1} (2-a)_k}\]
\[+\frac{p!}{1-a}\,\frac{\g(d)\g(d-a-c+p)}{\g(d-a)\g(d-c)(1-c)_p}\,\frac{\g(2-a)}{\g(1-a+p)}\]
upon use of the identity $(1-a)_{k+1}=(1-a) (2-a)_k$.

Now let $a=1-\epsilon$, where $\epsilon\to0$. Then
\[{}_3F_2\bl(\begin{array}{c}1-\epsilon,c,1\\d,p+1\end{array}\!\!;1\br)=\frac{p}{\epsilon}\sum_{k=0}^{p-1}(-)^k k!
\bl(\!\!\begin{array}{c}p-1\\k\end{array}\!\!\br)\,\frac{(1-d)_{k+1} \g(1+\epsilon)}{(1-c)_{k+1} \g(k+1+\epsilon)}\]
\[+\frac{p!}{\epsilon}\,\frac{\g(d) \g(d-c+p-1+\epsilon)}{\g(d-1+\epsilon)\g(d-c) (1-c)_p}\,\frac{\g(1+\epsilon)}{\g(p+\epsilon)}\]
\[=\frac{p}{\epsilon}\sum_{k=0}^{p-1}(-)^k 
\bl(\!\!\begin{array}{c}p-1\\k\end{array}\!\!\br)\,\frac{(1-d)_{k+1}}{(1-c)_{k+1}}\bl\{1-\epsilon(\gamma+\psi(k+1))+O(\epsilon^2)\br\}\]
\[+\frac{p}{\epsilon}\,\frac{\g(d)\g(d-c+p-1)}{\g(d-1)\g(d-c)(1-c)_p}\bl\{1-\epsilon(\gamma+\psi(p)-\psi(d-c+p-1)+\psi(d-1))+O(\epsilon^2)\br\}\]
upon use of the expansion $\g(z+\epsilon)=\g(z) \{1+\epsilon \psi(z)+O(\epsilon^2)\}$ and $\psi(1)=-\gamma$, where 
$\gamma=0.55721\ldots$ is the Euler-Mascheroni constant.

By means of the evaluation
\[\sum_{k=0}^{p-1}(-)^k\bl(\!\!\begin{array}{c}p-1\\k\end{array}\!\!\br)\,\frac{(1-d)_{k+1}}{(1-c)_{k+1}}=-\frac{\g(d)\g(d-c+p-1)}{\g(d-1)\g(d-c)(1-c)_p},\]
we then find that
\[{}_3F_2\bl(\begin{array}{c}1,1,c\\d,p+1\end{array}\!\!;1\br)=\frac{p\g(d)\g(d-c+p-1)}{\g(d-1)\g(d-c) (1-c)_p}\{\psi(d-c+p-1)-\psi(d-1)\}\]
\[+p\sum_{k=0}^{p-1}(-)^k\bl(\!\!\begin{array}{c}p-1\\k\end{array}\!\!\br)\,\frac{(1-d)_{k+1}}{(1-c)_{k+1}} \{\psi(p)-\psi(k+1)\}.\]
Now let $p\to n+1$ and replace the summation index $n$ in the finite sum by $n-k$, followed by a little rearrangement, to obtain 
\[{}_3F_2\bl(\begin{array}{c}1,1,c\\d,n+2\end{array}\!\!;1\br)=\frac{(n+1) \g(d)}{(1-c)_{n+1}}\bl\{\frac{(d-c)_n}{\g(d-1)}\{\psi(d-c+n)-\psi(d-1)\}\]
\bee\label{e31}
-\sum_{k=1}^n \bl(\!\!\begin{array}{c}n\\k\end{array}\!\!\br)\bl(\!\!\begin{array}{c}n\!+\!1\!-\!c\\k\end{array}\!\!\br)
\frac{k!\{\psi(n+1)-\psi(n+1-k)\}}{\g(d\!-\!n\!-\!1\!+\!k)}
\ee
for $n=0, 1, 2, \ldots\,$. This is the result stated in Theorem 1.\hfill $\Box$

\vspace{0.6cm}

\begin{center}
{\bf 3. \ Proof 2}
\end{center}
\setcounter{section}{3}
\setcounter{equation}{0}
\renewcommand{\theequation}{\arabic{section}.\arabic{equation}}
In a recent paper \cite{P}, the author examined convergent series expansions for the sum 
\[S(a,b)=\Lambda \sum_{m=1}^\infty \frac{J_\mu(am) J_\nu(bm)}{m^\al},\qquad \Lambda=\frac{2^{\mu+\nu}}{a^\mu b^\nu},
\]
where $J_\nu(x)$ is the Bessel function of the first kind and the factor $\Lambda$ is added for convenience. It is supposed that $\mu, \nu\geq0$, $a, b>0$  and $\al>0$ for absolute convergence. This was carried out using a Mellin transform approach and evaluation of residues.

In the particular case $\vartheta:=\al-\mu-\nu=2n+1$, $n=0, 1, 2 \ldots$ terms in $\log\,a$ appear in the expansions resulting from a double pole. When $a=b$, it is found that \cite[Theorem 2]{P}
\bee\label{e21}
S(a,a)=\mathop{\sum_{m=0}^\infty}_{\scriptstyle m\neq n}A_m (\fs a)^{2m}+\frac{(-)^n(\fs a)^{2n}\g(\al)\,\Upsilon_n(a)}{ \g(n\!+\!1\!+\!\mu)\g(n\!+\!1\!+\!\nu)\g(n\!+\!1\!+\!\mu\!+\!\nu)n!}
\end{equation}
for $0<a\leq\pi$, where
\[A_m=\frac{(-)^m}{m!}\,\frac{\g(1+\mu+\nu+2m)\,\zeta(\vartheta-2m)}{\g(1+\mu+m) \g(1+\nu+m) \g(1+\mu+\nu+m)},\]
and
\[\Upsilon_n(a)=\gamma-\log\,(\fs a)-\psi(\al)+\fs\psi(n+1)+\fs\psi(n+1+\mu)\hspace{3cm}\]
\[\hspace{6cm}+\fs\psi(n+1+\nu)+\fs\psi(n+1+\mu+\nu),\]
with $\zeta(s)$ being the Riemann zeta function.

When $a\neq b$ (with $a>b$), it is found that \cite[Theorem 4]{P}
\bee\label{e22}
S(a,b)=\frac{1}{\g(1+\nu)}\mathop{\sum_{m=0}^\infty}_{\scriptstyle m\neq n}B_m (\fs a)^{2m}+\frac{(-)^n (\fs a)^{2n}}{\g(1+\nu) \g(n+1+\mu) n!}\bl\{{\hat \Upsilon}_n(a) F_n(\mu,\chi)-\frac{1}{2}\Delta_n(\chi)\br\}
\ee
for $0<a+b\leq2\pi$, where $\chi:=b^2/a^2$. The coefficients $B_m$ are
\[B_m=\frac{(-)^m\zeta(\vartheta-2m)}{m! \g(1+\mu+m)}\,F_m(\mu,\chi),\qquad F_m(\mu,\chi)={}_2F_1\bl(\begin{array}{c}-m,-m-\mu\\1+\nu\end{array}\!\!;\chi\br),\]
where ${}_2F_1$ is the Gauss hypergeometric function. The other quantities appearing in (\ref{e22}) are
\[{\hat \Upsilon}_n(a)=\gamma-\log\,\fs a+\fs\psi(n+1+\mu)+\fs\psi(n+1)\]
and
\[
\Delta_n(\chi)=\sum_{k=1}^n \bl(\!\!\begin{array}{c}n\\k\end{array}\!\!\br) \bl(\!\!\begin{array}{c}n+\mu\\k\end{array}\!\!\br)\frac{{\cal D}_k(n,\mu) \,\chi^{k}}{(1+\nu)_k}\hspace{5cm}\]
\bee\label{e23}
\hspace{4cm}+\frac{(\mu)_{n+1} \chi^{n+1}}{(1+\nu)_{n+1} (n+1)}\,{}_3F_2\bl(\begin{array}{c}1,1,1-\mu\\n+\nu+2,n+2\end{array}\!\!;\chi\br),
\ee
where 
\[{\cal D}_k(n,\mu)=D_k(n)+k! \{\psi(n+1+\mu)-\psi(n+1+\mu-k)\}\]
with $D_k(n)$ defined in (\ref{e12}). 

If we now let $\chi=1$ and use the Gauss summation formula to show that
\[F_m(\mu,1)=\frac{\g(1+\nu) \g(1+\mu+\nu+2m)}{\g(1+\nu+m) \g(1+\mu+\nu+m)},\]
then it is easily seen that $A_m=B_m/\g(1+\nu)$ ($m\neq n$). Hence, from (\ref{e21}) and (\ref{e22}),
\[\frac{\g(\al) \Upsilon_n(a)}{\g(1+\nu_n)\g(1+\mu+\nu+n)}=\frac{1}{\g(1+\nu)}\bl\{{\hat\Upsilon}_n(a) F_n(\mu,1)-\frac{1}{2}\Delta_n(1)\br\}\]
from which we obtain
\[\Delta_n(1)=\frac{\g(\al)\g(1+\nu)}{\g(1+\nu+n)\g(1+\mu+\nu+n)} \bl\{2\psi(\al)-\psi(1+\nu+n)-\psi(1+\mu+\nu+n)\br\}.\]
From (\ref{e23}) we then find that
\[{}_3F_2\bl(\begin{array}{c}1,1,1-\mu\\n+\nu+2,n+2\end{array}\!\!;1\br)\hspace{8cm}\]
\[=\frac{(n+1) (1+\nu)_{n+1}}{(\mu)_{n+1}}\bl\{\frac{(1\!+\!\mu\!+\!\nu\!+\!n)_n}{(1+\nu)_n}\bl\{2\psi(\al)-\psi(1+\nu+n)-\psi(1+\mu+\nu+n)\br\}\]
\bee\label{e24}
-\sum_{k=1}^n \bl(\!\!\begin{array}{c}n\\k\end{array}\!\!\br)\bl(\!\!\begin{array}{c}n+\mu\\k\end{array}\!\!\br)
\,\frac{{\cal D}_k(n,\mu)}{(1+\nu)_k}\br\}.
\ee

If we now make the substitutions $c=1-\mu$ and $d=n+\nu+2$ in (\ref{e24}), so that $1+\mu+\nu+n=d-c$, $1+\nu=d-n-1$
and $\al=d-c+n$, we obtain after a little algebra
\[{}_3F_2\bl(\begin{array}{c}1,1,c\\d,n+2\end{array}\!\!;1\br)=\frac{(n+1) \g(d)}{(1-c)_{n+1}}\bl\{\frac{(d-c)_n}{\g(d-1)}\{2\psi(d-c+n)-\psi(d-1)-\psi(d-c)\}\]
\[-\sum_{k=1}^n\bl(\!\!\begin{array}{c} n\\k\end{array}\!\!\br)\bl(\!\!\begin{array}{c} n\!+\!1\!-\!c\\k\end{array}\!\!\br)\,\frac{{\cal D}_k(n,1-c)}{\g(d-n-1+k)}\br\}.\]
Use of the identity
\[\sum_{k=1}^n\bl(\!\!\begin{array}{c} n\\k\end{array}\!\!\br)\bl(\!\!\begin{array}{c} n\!+\!1\!-\!c\\k\end{array}\!\!\br)\frac{k!\{\psi(n\!+\!2\!-\!c)-\psi(n\!+\!2\!-\!c\!-\!k)\}}{\g(d-n-1+k)}=\frac{(d-c)_n}{\g(d-1)}\{\psi(d-c+n)-\psi(d-c)\}\]
to remove some of the $\psi$ functions,
we then obtain the result stated in (\ref{e11}).

The result in (\ref{e24}) was obtained on the assumption that $\mu, \nu\geq0$. But by the principle of analytic continuation the result in (\ref{e11}) will hold in the domain of the parameters $c$ and $d$ where both sides of (\ref{e11}) are analytic.
This concludes the proof. \hfill $\Box$
\bigskip

When $n=0$ and $n=1$ we have the summations
\[{}_3F_2\bl(\begin{array}{c}1,1,c\\d, n+2\end{array}\!\!;1\br)=\left\{\begin{array}{ll}
\dfrac{(d-1)}{1-c}\{\psi(d-c)-\psi(d-1)\} & (n=0)\\
\\
\dfrac{2(d-1)(d-c)}{(1-c)_2} \{\psi(d-c+1)-\psi(d-1)\} &\\
\hspace{6cm}+\dfrac{2(d-1))}{c-1} & (n=1),
\end{array}\right.\]
where we have used the property $\psi(z+1)=\psi(z)+1/z$. These summations agree with those presented in \cite[p.~452]{PBM}. In the case when $c=1$ ($n=0$) and $c=1, 2$ ($n=1$), the summations are given by their limiting forms. Numerical checks with the aid of {\it Mathematica} have confirmed the validity of (\ref{e11}).

\end{document}